\documentclass[12pt]{article}
\usepackage{latexsym}
\usepackage{amsmath}
\usepackage{amsfonts}
\usepackage{cite}

\newtheorem{Theorem} {Theorem} [section]
\newtheorem{Proposition} [Theorem] {Proposition}
\newtheorem{Lemma} [Theorem] {Lemma}
\newtheorem{Definition} [Theorem] {Definition}
\newcommand{\Proof}{\noindent{\bf Proof.}\quad }
\newcommand{\qed}{\hfill$\Box$\medskip}
\newcommand{\qedns}{{\hfill\hphantom{.}\nobreak\hfill$\Box$}}
\newcommand{\qedqed}{\hfill$\Box\Box$\medskip}

\newcommand{\Ff}{{\mathbb F}}

\newcommand{\adj}{\sim}
\newcommand{\PG}{{{\rm PG}}}
\newcommand{\<}{\langle}
\renewcommand{\>}{\rangle}
\newcommand{\linspan}[1]{{\langle{#1}\rangle}}
\newcommand\qnom[2]{\genfrac[]{0pt}{}{#1}{#2}}

\newcommand\mysqueeze[2]{%
\newdimen\origspacing%
\newdimen\newspacing%
\origspacing=\fontdimen2\font%
\setlength{\newspacing}{1\origspacing}%
\addtolength{\newspacing}{-#1}%
\fontdimen2\font=\newspacing%
{#2}%
\fontdimen2\font=\origspacing}

\title{The unique coclique extension property for apartments of buildings}
\author{A. E. Brouwer, J. Draisma \& \c{C}. G\"uven}
\date{16 November 2022}

\begin{document}
\maketitle

\begin{abstract}
We show that the Kneser graph of objects of a fixed type
in a building of spherical type has the unique coclique extension property
when the corresponding representation has minuscule weight and also
when the diagram is simply laced and the representation is adjoint.
\end{abstract}



\section{Introduction}
While studying generalizations of the Erd\H{o}s-Ko-Rado theorem
(\S\ref{sec:EKR}) and the chromatic number of Kneser-type graphs,
one needs information about maximal cocliques of near-maximum size,
cf.~\cite{bb,bbg,bbs,bbs2,dmw,cicekthesis,metsch2019,metsch2019a,mw}.
In this note we describe a simple construction that in the most
interesting cases produces all such near-maximum cocliques.

We say that a pair $(\Gamma,\Sigma)$ consisting of a graph and
an induced subgraph has the {\em unique coclique extension property}
when every maximal coclique $C$ of $\Sigma$ is contained in a unique
maximal coclique $D$ of $\Gamma$. Then $D$ necessarily consists of all
vertices $x$ of $\Gamma$ such that $C \cup \{x\}$ is a coclique, and
this property claims that $D$ thus defined is, indeed, a coclique.

As it turns out, Erd\H{o}s-Ko-Rado type results for set systems
and for their $q$-analogs (systems of subspaces) are very similar.
The unique coclique extension property explains part of this
similarity, giving a map from cocliques in an apartment of a building
to cocliques in the building.

Our Kneser graphs $\Gamma$ are defined on the objects (flags)
of some fixed type $J$ in a building of spherical type,
adjacent when they are in opposite chambers.
The subgraphs $\Sigma$ are those induced on such objects in
a fixed apartment.

Given a group $G$ with BN-pair $(B,N)$ (cf.~\cite{bourbaki})
and Weyl group $(W,R)$ and $J \subseteq R$ one has a building
with set of chambers $G/B$ and standard apartment $WB/B$.
If the building is spherical then
$W$ is finite and has a longest element $w_0$.
Let $J \subseteq R$ and put $X = \langle R\setminus J \rangle$,
a subgroup of $W$. Put $P = BXB$, a parabolic subgroup of $G$.
The vertices of the Kneser graph $\Gamma$ of type $J$ are the cosets $gP$
with $g \in G$, where $gP$ is adjacent to $hP$ when $Phg^{-1}P = Pw_0P$.
The vertices of $\Sigma$ are the cosets $wP$ with $w \in W$.

Let now $G$ be a reductive linear algebraic group (cf.~\cite{borel})
defined over the field $F$. Then the group $G(F)$ has a BN-pair.
Algebraic groups give rise to buildings with further structure,
and we can use their representation theory. We shall find that
the unique coclique extension property holds in cases
where the representation has minuscule weight or where the diagram
is simply laced and the representation is adjoint.


\begin{Theorem} \label{thm:Main}
Let $X_{n,i}$ be one of $A_{n,i}$ $(1 \le i \le n)$, $B_{n,1}$,
$B_{n,n}$, $C_{n,1}$, $D_{n,1}$, $D_{n,2}$, $D_{n,n}$, $E_{6,1}$,
$E_{6,2}$, $E_{7,1}$, $E_{7,7}$, $E_{8,8}$, $G_{2,1}$.

Let $\Gamma$ be the Kneser graph on the objects of type $i$ in a
building belonging to a Chevalley group with diagram $X_n$,
and let $\Sigma$ be the subgraph induced on an apartment.
Then the pair $(\Gamma,\Sigma)$ has the unique coclique extension property.
The same holds when $\Gamma$ is the Kneser graph on the objects
of type $\{1,n\}$ in a building with diagram $A_n$.
\end{Theorem}

\vskip -0.5em
The proof of this theorem is spread out over several sections.
In \S\ref{sec:Bilinear} we present a simple linear algebra argument
that will be used to establish most cases of our theorem, applicable
when adjacency in $\Gamma$ corresponds to the non-vanishing of
a certain bilinear map.
The cases $A_{n,i}$, $A_{n,\{1,n\}}$, $B_{n,1}$/$C_{n,1}$/$D_{n,1}$,
$D_{n,2}$ are treated in
\S\S\ref{sec:Ani}, \ref{sec:PH}, \ref{sec:PointPolar}, \ref{sec:LineOrtho},
respectively, using geometric arguments. Using representation theory of
algebraic groups, the cases $B_{n,n}$, $D_{n,n}$, $E_{6,1}$, and $E_{7,7}$
(corresponding to a representation with minuscule weight)
follow in \S\ref{sec:Minuscule}, while the cases $D_{n,2}$, $E_{6,2}$,
$E_{7,1}$, $E_{8,8}$, and $A_{n,\{1,n\}}$ (where the representation is adjoint)
follow in \S\ref{sec:Adjoint}.

Some examples of cases where the unique coclique extension property fails
are given in \S\ref{Nonexamples}. Given some examples or nonexamples,
\S\ref{Varying} discusses further derived (non)examples. For example,
the unique coclique extension property holds for type $A_{n,\{1,2\}}$
($n \ge 2$).


\section{The Erd\H{o}s-Ko-Rado theorem}\label{sec:EKR}
Let $n,k$ be positive integers with $n > 2k$. The Erd\H{o}s-Ko-Rado
theorem \cite{EKR} affirms that an intersecting family of $k$-subsets
of an $n$-set has size at most $\binom{n-1}{k-1}$ (and specifies
the families attaining this bound).
The Hilton-Milner theorem \cite{hiltonmilner}
specifies the size (namely $\binom{n-1}{k-1}-\binom{n-k-1}{k-1}+1$)
and structure of the second largest maximal intersecting families.

The Kneser graph $K(n,k)$ is the graph on the $k$-subsets of an $n$-set,
adjacent when disjoint. These two theorems give the size and structure
of the largest and next largest maximal cocliques in the Kneser graph.

These theorems have $q$-analogs.  Hsieh \cite{hsieh} showed that a system
of nontrivially intersecting $k$-subspaces of an $n$-space over $\Ff_q$
has size at most $\qnom{n-1}{k-1}$.
and in \cite{BBCFMPS} it was shown that the next largest size
is $\qnom{n-1}{k-1}-q^{k(k-1)}\qnom{n-k-1}{k-1}+q^k$.

\medskip
This setup can be generalized. Given a spherical building (that is, one
with a finite Coxeter group $W$), the Kneser graph $\Gamma$ of type $J$
has as vertices the objects of type $J$, adjacent when they are
in opposite chambers.\footnote{
Two singular $k$-subspaces $A,B$ of a polar space of rank $n$ are adjacent
in the Kneser graph when $A^\perp \cap B = 0$, a relation that is stricter
than the disjointness $A \cap B = 0$.
}
Parameter information for these graphs was given in \cite{qpows,dmm}.

The Kneser graph $K(n,k)$ and its $q$-analog are the special cases
of the Kneser graph of type $k$ in an apartment or building of type $A_{n-1}$.

\medskip
For analogs for maximal singular subspaces in a polar space of rank $n$,
see \cite{Stanton1980}, \cite{pepe-et-al}.
Tanaka \cite{Tanaka} classifies the subsets in certain association schemes
of which width plus dual width equals the number of classes, and derives
EKR-type results.
For a book-length discussion of EKR-type results, see \cite{GodsilMeagher}.

\medskip
In the cases where the unique coclique extension property holds
it often yields the largest cocliques in the Kneser graph.
For example, for $A_{n,\{1,n\}}$ (over $\Ff_q$) the maximum size is
$1+2q+3q^2+\cdots+nq^{n-1}$ and examples of this size
arise using this construction (see \S\ref{sec:PH}).

\section{The coclique extension lemma} \label{sec:Bilinear}

Let $V,V'$ be vector spaces and $\mu : V \times V \to V'$
a bilinear map that is reflexive (i.e., satisfies $\mu(x,y) = 0$
if and only if $\mu(y,x) = 0$). For $A \subseteq V$, write $A^\perp =
\{v \in V \mid \mu(v,a) = 0 ~\mbox{\rm for all}~a \in A\}$.

Suppose the vertex set $V\Gamma$ of a graph $\Gamma$ can be embedded
into the projective space $PV$ corresponding to $V$ in such a way that
vertices $x,y$ are nonadjacent in $\Gamma$ if and only if $\mu(x,y) = 0$
(and in particular $\mu(x,x) = 0$ for all vertices).

\begin{Lemma}
Let $C$ be a coclique in $\Gamma$, so that $C \subseteq C^\perp$,
and let $D$ be the set of vertices $x$ of $\Gamma$ for which
$C \cup \{x\}$ is also a coclique, so that $D = C^\perp \cap V\Gamma$.
If $D$ is contained in $\langle C \rangle$, the linear span of $C$ in $V$,
then $D$ is a coclique. \qedns
\end{Lemma}

This will be our main tool in many of the cases.

\section{Subspaces of a projective space}\label{sec:Ani}
\label{sec:Projective}

Fix integers $m,n$, where $1 \le m \le \frac{1}{2}n$.
Let $N = \{1,\ldots,n\}$.
Let $F$ be a division ring, let $V$ be an $n$-dimensional
left vector space over $F$ with basis $\{e_1,\ldots,e_n\}$, and let
$PV$ be the corresponding projective space.
Each $m$-set $K \subseteq N$ determines an $m$-space $\phi K \le PV$
by $\phi K = \langle e_i \mid i \in K \rangle$.

\begin{Proposition}\label{An}
Let $C$ be a maximal collection of pairwise intersecting
$m$-subsets of $N$.
Let $D$ be the collection of all $m$-spaces in $PV$ that meet
$\phi K$ for all $K \in C$. Then $D$ is a maximal collection
of pairwise intersecting $m$-subspaces of $PV$.
\end{Proposition}
\vspace{-2mm}
(Here and elsewhere, dimensions are vector space dimensions.)

\medskip\noindent
{\bf First proof} (for the case of a field $F$). 
Let $\bigwedge V$ be the exterior algebra of $V$.
Map $m$-subspaces of $PV$ to projective points in
$P(\bigwedge V)$ via
$$\psi : U = \langle u_1,\ldots,u_m \rangle \mapsto
\langle u_1 \wedge \ldots \wedge u_m \rangle .$$
Now $U \cap U' \ne 0$ if and only if $\psi U \wedge \psi U' = 0$.
For $M = \{i_1,\ldots,i_m\} \subseteq N$ with $i_1 < \ldots < i_m$
let $e_M = e_{i_1} \wedge \ldots \wedge e_{i_m}$.
The $e_M$ form a basis for the degree $m$ part of $\bigwedge V$.
Suppose $U \in D$, where $\psi U = \langle \,\sum \alpha_M e_M \,\rangle$.
Then $\alpha_M = 0$ whenever $M$ is disjoint from some $K
\in C$. Indeed, if $K \cap M = \emptyset$, then the coefficient of
$e_{K \cup M}$ in $\psi U \wedge \psi\phi K$ is $\alpha_M$
up to a nonzero constant. 

So, if $\alpha_M \ne 0$ then $M$ meets all $K \in C$, and,
since $C$ is maximal, $M \in C$.
We see that the 1-space $\psi U$ is contained in $\langle \psi\phi C \rangle$.
If $U,U' \in D$, then $\psi U, \psi U'\subseteq \langle \psi\phi C \rangle$
implies $\psi U \wedge \psi U' = 0$ so that $U \cap U' \ne 0$.
Maximality of $D$ is clear.
\qed

Note that the proof above implements the argument from
\S\ref{sec:Bilinear}. If $F$ is not commutative, we need a different argument. 

\medskip\noindent
{\bf Second proof} (for the general case).
Given a matrix $A$ with entries in $F$, let its {\em rank} be
the dimension of the left vector space spanned by its rows,
or, equivalently, the dimension of the right vector space
spanned by its columns.

Given a matrix $A$ with entries in $F$ and columns indexed by $N$,
we find a matroid $M_A$ on $N$ by letting the rank of $M \subseteq N$
be the rank of the submatrix of $A$ on the columns indexed by $M$.

Since $V$ has a fixed basis $(e_i)_i$, we can regard each $v \in V$ as
a row vector. Accordingly, we regard $V$ as a left vector space over
$F$. Given a (left $F$-) subspace $U$ of $V$, let $A(U)$ be a matrix
of which the rows form a basis of $U$. We find a matroid $M_U$ on $N$
by taking $M_U = M_{A(U)}$. It is independent of the choice of $A(U)$.

An $m$-space $U$ meets $\phi K$ if and only if no basis of $M_U$ is
disjoint from $K$. Hence $U \in D$ if and only if each basis of $M_U$
meets all elements of $C$, i.e., by maximality of $C$,
if and only if each basis of $M_U$ belongs to $C$.
Now the claim follows from:

\medskip\noindent
{\bf Lemma}. {\em Let $U,U'$ be disjoint subspaces of $V$.
Then $M_U, M_{U'}$ have disjoint bases.}

\noindent
{\bf Proof}.
Let $M$ be the union matroid of $M_U,M_{U'}$, that is the matroid
of which the independent sets are the unions of an independent set
of $M_U$ and one of $M_{U'}$. As is well known
(see, e.g., \cite{Lex}, Ch.~42), the rank function of $M$ is given by
$$r(K) = \min_{L \subseteq K} \,(|K \setminus L| + r_U(L) + r_{U'}(L)) ,$$
where $r_U,r_{U'}$ are the rank functions of $M_U,M_{U'}$.
Let $\dim U = m$, $\dim U' = m'$. We have to show that
$r(N) = m+m'$, that is, that
$$n - |L| + r_U(L) + r_{U'}(L) \ge m+m'$$ for all $L \subseteq N$.
But $r_U(L) + r_{U'}(L) \ge r_{U+U'}(L)$ and
$m+m' - r_{U+U'}(L) \le n - |L|$ since the subspace of $U+U'$
consisting of the vectors vanishing on the positions in $L$
is contained in the subspace of $V$ of such vectors.
\qedqed

\section{{\hskip -1.7pt}Point-hyperplane flags in a
projective space}\label{sec:PH}
In this section we study the Kneser graphs for type $A_{n,\{1,n\}}$.

Let $V$ be an $(n+1)$-dimensional left vector space over the
division ring $F$, and let $\Gamma$ be the graph of which
the vertices are the point-hyperplane flags $(P,H)$
(where points are 1-spaces), where $(P,H) \adj (Q,I)$
when $P \not\subseteq I$ and $Q \not\subseteq H$.

Let $V$ have basis $\{ e_i \mid i \in N \}$, where $N = \{0,\ldots,n\}$.
Let $\Sigma$ be the subgraph of $\Gamma$ induced on the set of $n(n+1)$
vertices $\{ (i,N\setminus j) \mid i \ne j \}$, where $(i,N\setminus j)$
abbreviates $(\linspan{e_i},\linspan{\{e_h \mid h \ne j\}})$.
A maximal coclique $C$ in $\Sigma$ has size $\binom{n+1}{2}$ 
and contains precisely one vertex from every pair
$(i,N\setminus j)$, $(j,N\setminus i)$.

Fix $C$, and let $D$ be the set of vertices of $\Gamma$ not adjacent to
any vertex in $C$. We claim that $D$ is a coclique in $\Gamma$.
Indeed, suppose that $D$ contains adjacent vertices
$(\linspan{p},I)$ and $(\linspan{q},J)$. Let
$p = \sum \alpha_i e_i$ and $q = \sum \beta_j e_j$.
Since $p \notin J$ there is an $i$ with $\alpha_i \ne 0$
and $e_i \notin J$. Similarly, there is a $j$ with $\beta_j \ne 0$
and $e_j \notin I$.
Since $p \in I$ we have $i \ne j$.
W.l.o.g., let $(i,N\setminus j) \in C$.
Both $e_i \in J$ and $q \in \linspan{N\setminus j}$ are impossible,
contradiction.
It follows that the unique coclique extension property holds here.

Now let $F = \Ff_q$. It was shown in \cite{bbg} that the largest cocliques
in $\Gamma$ have size $1+2q+3q^2+ \cdots + n q^{n-1}$,
that if we pick $C$ so that it contains $(i,N\setminus j)$ precisely
when $i < j$, then $D$ will be a maximum size coclique in $\Gamma$,
and that all maximum cocliques in $\Gamma$ are obtained in this way
(for a suitable choice of the apartment $\Sigma$ and ordering of $N$).

\section{Points in a polar space} \label{sec:PointPolar}
%
A {\em polar space of rank $n$} is a building of type $C_n$
with points and lines being the objects of types 1 and 2.
All objects are singular subspaces.
The Kneser graph $\Gamma$ on the points is the noncollinearity graph.
An apartment in this building is a cross polytope on $2n$ vertices
(for $n=2$ a square, for $n=3$ an octahedron),
and the Kneser graph $\Sigma$ is $nK_2$, the disjoint union
of $n$ copies of $K_2$.

Every maximal coclique $C$ of $\Sigma$ has size $n$, and is
contained in a unique maximal singular subspace $D$,
that is, in a unique maximal coclique of $\Gamma$.
It follows that the unique coclique extension property holds.

The buildings of Chevalley or twisted Chevalley groups of types
$B_n$, $C_n$, $D_n$, ${}^2A_{2n-1}$, ${}^2A_{2n}$, and ${}^2D_{n+1}$
are polar spaces of rank $n$.

\medskip\noindent
For the classical generalized hexagon of type $G_2$, the points
are the points of the polar space of type $B_3$, and the Kneser graphs
of both coincide, and so do the subgraphs induced on an apartment.
It follows that the unique coclique extension property also holds
for type $G_{2,1}$.

\section{{\hskip -2.8pt}Totally singular lines in an
orthogonal space} \label{sec:LineOrtho}
In this section we study the Kneser graphs for type
$D_{n,2}$.

Fix an integer $n \ge 2$. Let $N = \{1,1',2,2',\ldots,n,n'\}$.
It is provided with an involution $'$ that interchanges $i$ and $i'$.
Let $V$ be a $2n$-dimensional vector space (over some field $F$)
with basis $\{e_s \mid s \in N \}$ and let $Q$ be the quadratic form on $V$
defined by $Q(x) = x_1x_{1'} + \ldots + x_nx_{n'}$. 
Now $(V,Q)$ is a nondegenerate orthogonal space with maximal Witt index.
Each pair $\{s,t\} \subseteq N$ determines a line in $PV$ by
$\phi \{s,t\} = \langle e_s, e_t \rangle$.
This line will be totally singular when $t \ne s'$.


\begin{Proposition}
Let $C$ be a maximal collection of pairs $\{s,t\}$ in $N$ such that
$t \ne s'$ and if $\{s,t\} \in C$ then $\{s',t'\} \not\in C$.
Let $D$ be the collection of all totally singular lines in $PV$
that meet $(\phi K)^\perp$ for all $K \in C$.
Then $L^\perp \cap M \ne 0$ for any two lines $L,M \in D$.
\end{Proposition}
\vspace{-1mm}
Note that $L^\perp \cap M \ne 0$ if and only if $L \cap M^\perp \ne 0$.

\medskip\Proof
As before, let $\psi$ map subspaces of $PV$ to points in $P(\bigwedge V)$.
For two lines $L,M$ we have $L \cap M^\perp \ne 0$ whenever
$\psi L \wedge \psi (M^\perp) = 0$.
For $K \in C$ we have $\psi ((\phi K)^\perp) = e_{N \setminus K'}$.
It follows that the line $L$ with $\psi L = \sum \alpha_P e_P$
(with $P$ running over the 2-subsets of $N$) intersects
$(\phi K)^\perp$ for $K \in C$ if and only if $\alpha_{K'} = 0$.
Let $L \in D$, and suppose $\psi L$ involves $e_{i,i'}$
for some $i$. Since $L$ is totally singular, it follows that $\psi L$
also involves $e_{j,j'}$ for some $j \ne i$. The collection $C$ contains
precisely one of $\{s,t\}$ and $\{s',t'\}$ whenever $t \ne s'$, so we may
assume that $C$ contains $\{i,j\}$ and $\{i,j'\}$, so that
$\alpha_{\{i',j\}} = \alpha_{\{i',j'\}} = 0$. Let $L = \langle u,v \rangle$
with $u = \sum u_s e_s$ and $v = \sum v_s e_s$ where $u_{i'} = 0$.
Then $v_{i'} \ne 0$ and $u_j = u_{j'} = 0$ so that $\psi L$ does not
involve $e_{j,j'}$, a contradiction. So, $\alpha_P = 0$ when $P = \{i,i'\}$.
We proved that $\psi L \in \langle \psi\phi C\rangle$.
Now the proof finishes as the first proof of Proposition \ref{An}.
\qed

The above proved the unique coclique extension property for the Kneser
graph of type $D_{n,2}$. It also holds for the disjointness graph in
this geometry.

\begin{Proposition}
Let $C$ be a maximal collection of pairwise intersecting $2$-subsets
$\{s,t\}$ of $N$ with $t \ne s'$.
Let $D$ be the collection of all totally singular lines in $PV$
that meet $\phi K$ for all $K \in C$. Then $D$ is a maximal collection
of pairwise intersecting totally singular lines in $PV$.
\end{Proposition}
\Proof
This is trivial: one finds either all lines in a totally singular plane
or all totally singular lines on a fixed point.
\qed

\section{Minuscule weights} \label{sec:Minuscule}
Let $G$ be a split connected reductive algebraic group defined
over a field $F$, $T$ a maximal split torus, and $V$ an irreducible
algebraic representation of $G$ over $F$.
For each character $\chi$ of $T$, let $V_\chi$ be the subspace of
$V$ on which $T$ acts with that character. Assume that $V$ is minuscule,
i.e., that $W$ acts transitively on the weights $\chi$ for which $V_\chi$
is nonzero. Then it follows from irreducibility that all weight spaces
are $1$-dimensional. Fix such a weight $\lambda$, let $v_\lambda \in
V_\lambda(F)$ be nonzero, and let $P$ be the stabilizer of $v_\lambda$.

Now the Kneser graph $\Gamma$ is the graph with vertex set $G(F)/P(F)$,
embedded in the projective space corresponding to $V$ as the orbit of
$\langle v_\lambda \rangle$ under $G(F)$, and $V\Sigma$ is the orbit
of $\langle v_\lambda \rangle$ under the Weyl group $W$. The neighbours
of a vertex in $V\Sigma$, corresponding to the weight $\lambda' \in W$,
are the vertices corresponding to weights that are at maximal distance
from $\lambda'$.

We will look for reflexive bilinear maps as in the following definition.

\begin{Definition}
A reflexive bilinear map $V \times V \to V'$, where $V,V'$ are
$G$-modules, is called {\em good} if it has the following
properties: two vertices $x,y$ of $\Gamma$ are
adjacent in $\Gamma$ if and only if $\mu(x,y) \ne 0$; and, 
moreover, for each $x \in V\Sigma$, the nonzero vectors $\mu(x,y)$
with $y \in V\Sigma$ are linearly independent. 
\end{Definition}

\begin{Proposition}
(i) For each minuscule representation $V$, there exists a good bilinear map
$\mu:V \times V \to V'$. 

(ii) The pair $(\Gamma,\Sigma)$ has the unique coclique extension property.
\end{Proposition}
\vspace{-2mm}
This settles the cases $A_{n,j}$ $(1 \le j \le n)$, $B_{n,n}$, $C_{n,1}$,
$D_{n,1}$, $D_{n,n}$, $E_{6,1}$, $E_{7,7}$.
\vspace{2mm}

\Proof
Let us do part (ii) first, given part (i).

Assume that we have a good bilinear map $\mu:V \times V \to V'$.
The module $V$ is a direct sum of 1-dimensional weight spaces permuted by
$W$ and spanned by the vertices of $\Sigma$, so that $\dim V = |V\Sigma|$.
Pick a basis vector $e_s$ for each $s \in \Sigma$.  For $x \in V\Gamma$
write $x = \sum x_s e_s$.  If $x \in C^\perp$, then $0 = \mu(x,e_c) =
\sum x_s \mu (e_s,e_c)$ for each $c \in C$. Since, by assumption, the
nonzero vectors $\mu(e_s,e_c)$ are linearly independent, we find that
for all $s$ with $x_s \neq 0$, we have $\mu(e_s,e_c)=0$ for all $c$.
Since $C$ is maximal, any $s$ with this property lies in $C$, and we
showed that $x \in \langle C \rangle$.

It remains to find a good bilinear map $\mu:V \times V \to V'$.

%

We already did $A_{n,j}$ (using the bilinear map
$\mu : V \times V \to \bigwedge V_0$, where $V = \bigwedge^j V_0$
and $V_0$ is the natural module for ${\rm SL}_{n+1}$ of dimension $n+1$),
and $C_{n,1}$ and $D_{n,1}$ (using the natural symplectic or
symmetric bilinear form $\mu : V \times V \to F$).

More generally, if $V$ is a self-dual representation, then there is
a $G$-equivariant reflexive bilinear form $\mu:V \times V \to F$. In
that case, every vertex in $V \Sigma$ has a unique neighbour in $V
\Sigma$, namely, the one with the opposite weight, so the condition on
linear independence is automatically satisfied. That settles $B_{n,n}$,
$D_{n,n}$ ($n$ even), and $E_{7,7}$, in addition to 
$C_{n,1}$ and $D_{n,1}$.

For $D_{n,n}$ ($n$ odd), $\Gamma$ is the graph on maximal isotropic
subspaces in a $2n$-dimensional orthogonal space $V_0$ that intersect one
given such space in an odd-dimensional subspace. Two are adjacent if they
intersect in a $1$-dimensional space.  Here $V$ is the spin representation
of the spin group $G$. The spin representation has a nonzero reflexive
$G$-equivariant bilinear map $\mu: V \times V \to V_0$ (symmetric when
$n$ is $1$ modulo $4$, skew-symmetric when $n$ is $3$ modulo $4$).%
\footnote{
If the characteristic is not 2, $\mu$ can be taken to be the map $\beta_1$
from Chevalley \cite{chevalley}, III.4.4, and a similar construction works
in characteristic 2.
Much more generally, the existence of $\mu$ follows from
\cite[Theorem 4.3.2]{brion-et-al}, as pointed out to us by Shrawan Kumar.
Our $V$ is the Weyl module $V(\omega_n)$, where the $n$-th fundamental weight
$\omega_n$ is the minuscule weight under consideration, and the theorem
there says that there is a unique surjective equivariant bilinear map
from the submodule $M$ of $V \otimes V$ generated by the tensor product
$v_{\omega_n} \otimes v_{w_0 \omega_n}$ into
$V(\overline{\omega_n + w_0 \omega_n})$. Here
$w_0$ is the longest element in the Weyl group, $v_\lambda$ spans the
(one-dimensional) weight space in $V$ corresponding to $\lambda$, and
$\overline{\lambda}$ is the unique dominant weight in the Weyl group orbit
of $\lambda$.  In the current case, a straightforward check shows that $M$
is all of $V \otimes V$, $w_0 \omega_n=-\omega_{n-1}$ (for $D_n$ with $n$
odd, $w_0$ is the order-two automorphism of the diagram followed by $-1$),
and $\overline{\omega_n - \omega_{n-1}}=\omega_1$, the highest weight
of the standard representation.
}
%
If $a,b \in V$ represent isotropic $n$-dimensional spaces $A,B \subseteq
V_0$ that intersect in an odd-dimensional space $C$, and if $\mu(a,b)$
is nonzero, then equivariance implies that $C$ is one-dimensional and
spanned by $\mu(a,b)$. Furthermore, if $e_1,e_{1'},\ldots,e_n,e_{n'}$
is the basis of $V_0$ discussed in \S\ref{sec:LineOrtho}, and $a \in V$
is the vector representing $\langle e_1,\ldots,e_n \rangle$, then the
$2^{n-1}$ basis vectors of the spin representation each represent a space
spanned by $e_i$ for an odd number of $i \in \{1,\ldots,n\}$, and $e_{i'}$
for the remaining $i$. It then follows that precisely $n$ of these
intersect $\langle e_1,\ldots,e_n \rangle$ in a $1$-dimensional space,
and as $b$ ranges over the corresponding basis vectors, $\mu(a,b)$ ranges
over $e_1,\ldots,e_n$. So these are linearly independent as desired.

For $E_{6,1}$, the module $V$ is 27-dimensional, and the dual module
$V'$ is that of $E_{6,6}$. There is a nonzero symmetric $G$-equivariant
bilinear map $V \times V \to V'$.%
\footnote{
Also this is a special case of \cite[Theorem 4.3.2]{brion-et-al}.
}
It has the property that $\mu(x,y) = z$, where $z$ is the unique
symplecton containing $x,y$ when $x,y$ are noncollinear,
and $z = 0$ otherwise.  The apartment is the Schl\"afli
graph (27 vertices, valency 16).  For fixed $c$ there are 10 nonadjacent
vertices, and the 10 symplecta $\mu(c,s)$ are distinct and hence linearly
independent.
\qed


\section{Adjoint representation} \label{sec:Adjoint}
Assume that the algebraic group $G$ has simply laced Dynkin diagram,
and consider its adjoint representation $V$, the Lie algebra of $G$.
In this case, $\Gamma$ is the orbit of root vectors, and two are adjacent
if they generate a Lie subalgebra isomorphic to $\mathfrak{sl}_2(F)$, i.e.,
with respect to a suitable choice of maximal torus, they correspond to
opposite roots.

\begin{Proposition}
There is a good bilinear form $V \times V \to F$,
and the unique coclique extension property holds. 
\end{Proposition}
\vspace{-1mm}
This covers $D_{n,2}$ $(n \ge 4)$, $E_{6,2}$, $E_{7,1}$, $E_{8,8}$,
and $A_{n,\{1,n\}}$ $(n \ge 2)$.

\medskip\noindent\Proof
The module $V = V(\lambda)$ carries the structure of a Lie algebra
${\mathfrak g}$. Let $\Phi$ be the root system. Then
$\dim V = |\Phi| + \ell$ where $\ell$ is the Lie rank of $G$,
the dimension of a Cartan subalgebra ${\mathfrak h}$.
Pick a Chevalley basis consisting of $e_s$ for $s \in \Phi$
and $e_i$ for $1 \le i \le \ell$, where $V\Sigma$ is the
collection of root spaces $\langle e_s \rangle$ for $s \in \Phi$,
and ${\mathfrak h}$ is spanned by the $e_i$.

\medskip
Let $r \in \Phi$, and consider $[e_r,[e_r,v]]$ for arbitrary $v \in V$.
If $v = e_s$ for some $s \in \Phi$, then $r,s$ span a 2-dimensional
root system, and we see that this is zero, unless $r+s = 0$, in which
case it is $-2e_r$. If $v = h \in {\mathfrak h}$, then
$[e_r,[e_r,v]] = 0$. More generally, for $v = h + \sum v_s e_s$,
we find $[e_r,[e_r,v]] = -2v_{-r}e_r$.
Let us call an element $x$ of a Lie algebra {\em extremal} when
$[x,[x,y]]$ is a multiple of $x$ for all $y$.
(The usual definition adds a requirement for the case of characteristic 2
that we don't need.)
We just observed that $e_r$ is extremal. It follows that $x$ is extremal
for every $x \in V\Gamma$.

\medskip
The graph $\Sigma$ has valency 1, so that a maximal clique $C$ contains
exactly one from each pair of opposite root vectors.

The usual Killing form has a constant factor that makes it identically zero
in small characteristics. Let $\mu$ be the reduced Killing form
(cf.~\cite{garibaldi-et-al}, \S8 and \cite{GrossNebe2004}, \S5),
then $\mu(e_r,e_{-r}) = 1$, so that the bilinear form $\mu$ is good.

\medskip
We want to show that if $x,y \in C^\perp \cap V\Gamma$, then $\mu(x,y) = 0$.
Let $x = x_0 + \sum x_s e_s$ and $y = y_0 + \sum y_s e_s$ with
$x_0,y_0 \in {\mathfrak h}$.
If $\< e_c \> \in C$, then $x_{-c} = \mu(x,e_c) = 0$, and similarly $y_{-c}=0$.
It follows that $\mu(x,y) = \mu(x_0,y_0)$.
If $x_0$ lies in the radical of $\mu$, then $\mu(x_0,y_0) = 0$.
%
%
If not, then there is a root $t \in \Phi$ with $t(x_0) \ne 0$.
If necessary replace $t$ by $-t$ to make sure that
$\langle e_t \rangle \not\in C$.
Since $x$ is extremal, $[x,[x,e_t]]$ is a multiple of $x$.
On the other hand,
the coefficient of $e_t$ in this expression is $t(x_0)^2 \ne 0$
(since $x_s x_{-s} = 0$ for all $s \in \Phi$).
This is a contradiction.
\qed

\section{Nonexamples}\label{Nonexamples} 
The unique extension property does not hold for
arbitrary types. We give counterexamples for types
$A_{n-1,\{i,n-i\}}$ $(1 < i < n/2)$,
$B_{3,2}$, $C_{3,3}$, $D_{4,\{3,4\}}$, and $D_{5,3}$.

\begin{Proposition}
The unique coclique extension property does not hold for
the Kneser graph on the flags of type $\{i,n-i\}$ in a building
with diagram $A_{n-1}$, where $1 < i < n/2$.
\end{Proposition}
\Proof
Let $V$ be a vector space of dimension $n$ with basis
$e_1,\ldots,e_n$. Put $u = e_1+e_2$, $v = e_1+e_n$.
\mysqueeze{0.8pt}{Let $A = \langle u,e_3,\ldots,e_{i+1} \rangle$
and $A' = \langle v,e_{n-1},\ldots,e_{n-i+1} \rangle$} 
so that $A$ and $A'$ are $i$-spaces.
Let $B = \langle u,e_3,\ldots,e_{n-i},e_n \rangle$
and $B' = \langle v,e_{n-1},\ldots,e_{i+2},e_2 \rangle$
so that $B$ and $B'$ are $(n-i)$-spaces.
Now the flags $F = (A,B)$ and $F' = (A',B')$ are adjacent since
$A \cap B' = A' \cap B = 0$.
Note that $F$ is mapped to $F'$ by the coordinate permutation
$(2,n)(3,n-1)\ldots$.

The graph $\Sigma$ has valency 1 and each choice of one vertex
from each edge of $\Sigma$ yields a maximal coclique $C$.
We can find a maximal coclique $C$ compatible with $F,F'$
when in no edge of $\Sigma$ both endpoints are adjacent to
either $F$ or $F'$.
Let $N = \{1,\ldots,n\}$. The vertices of $\Sigma$ are pairs
$(S_I,S_J)$ where $S_I = \langle e_i \mid i \in I \rangle$
and $|I| = i$, $|J| = n-i$, $I \subseteq J$.
The unique neighbour of $(S_I,S_J)$ is
$(S_{N\setminus J},S_{N\setminus I})$.
If $(S_I,S_J)$ is adjacent to $F$, then $S_I \cap B = S_J \cap A = 0$,
so that $I = \{j,n-i+1,\ldots,n-1\}$ and $J = \{j,i+2,\ldots,n\}$,
where $j \in \{1,2\}$.
If $(S_I,S_J)$ is adjacent to $F'$, then
$I = \{k,i+1,\ldots,3\}$ and $J = \{k,n-i,\ldots,2\}$
where $k \in \{1,n\}$. Altogether 4 vertices in $\Sigma$
are adjacent to either $F$ or $F'$, and for $i > 1$
this set of 4 does not contain any edge. \qed

\medskip\noindent
$B_{3,2}$ nonexample:
Let ${\rm char}~F \ne 2$ and let $V_0 = F^7$, with basis vectors
$e_i$, $1 \le i \le 7$.
Provide $V_0$ with the nondegenerate quadratic form
$Q(x) = x_1x_2 + x_3x_4 + x_5x_6 - x_7^2$.
The geometry of totally singular subspaces of $(V_0,Q)$ has type $B_3$.
For type $B_3$, the adjoint representation is $\wedge^2 V_0$,
corresponding to $B_{3,2}$, the geometry of totally singular lines.
Consider the apartment with points $\langle e_i \rangle$, $1 \le i \le 6$,
and let $C$ consist of the six totally singular lines
$13,14,25,26,35,36$, where $ij$ abbreviates $\langle e_i,e_j \rangle$.
Pick $u = e_1$, $v = e_3+e_4+e_7$ and
$u' = e_2$, $v' = e_5+e_6+e_7$. Then $\langle u,v \rangle$
and $\langle u',v' \rangle$ are totally singular and adjacent
in the Kneser graph and both nonadjacent to all vertices of $C$.

\medskip\noindent
$C_{3,3}$ nonexample:
Let ${\rm char}~F \ne 2$ and let $V_0 = F^6$, with basis vectors
$e_i$, $1 \le i \le 6$.
Provide $V_0$ with the nondegenerate symplectic form
$f(x,y) = x_1y_2-x_2y_1+x_3y_4-x_4y_3+x_5y_6-x_6y_5$.
The geometry of totally isotropic subspaces~of~$(V_0,f)$ has type $C_3$.
Consider the apartment with points $\langle e_i \rangle$, $1 \le i \le 6$,
and let $C$ consist of the four totally isotropic planes
$135,146,236,245$ where $ijk$ abbreviates $\langle e_i,e_j,e_k \rangle$.
The two planes $\langle e_1, e_3+e_6, e_4+e_5 \rangle$ and
$\langle e_1+e_3-e_5, e_2+e_6, e_4+e_6 \rangle$ are totally isotropic
and disjoint (since ${\rm char}~F \ne 2$), i.e., adjacent in the Kneser graph,
while nonadjacent to all vertices of $C$.

\medskip\noindent
Planes in $D_4$ nonexample:
Let $N = \{1,2,3,4,1',2',3',4'\}$ and let $'$ be the involution on $N$
that maps $i$ to $i'$.
Let ${\rm char}~F = 2$ and let $V_0 = F^8$, with basis vectors $e_i$
$(i \in N)$ and provided with the nondegenerate quadratic form
$Q(x) = x_1x_{1'} + x_2x_{2'} + x_3x_{3'} + x_4x_{4'}$.
The geometry of totally singular subspaces of $(V_0,Q)$ has type $D_4$.
Consider the apartment with points $\langle e_i \rangle$, $i \in N$.
Let $\pi = \langle e_1+e_2, e_{1'}+e_{2'}, e_4 \rangle$ and
$\pi' = \langle e_{1'}+e_3, e_1+e_{3'}, e_{4'} \rangle$.
Then $\pi,\pi'$ are totally singular planes, adjacent in the Kneser graph
on the planes in $D_4$ (since $\pi^\perp \cap \pi' = 0$).
Here $\Sigma$ has valency 1,
and one checks that there is no edge in $\Sigma$ such that both ends
are adjacent to $\pi$ or $\pi'$, so there is a maximal coclique $C$
compatible with both $\pi$ and $\pi'$,
and the unique coclique extension property fails.
This was the case $D_{4,\{3,4\}}$.
Of course this also means that it fails for $D_{n,3}$ for all $n \ge 5$.

\section{Varying $J$}\label{Varying}
Given a building of spherical type, let $\Gamma$, $\Gamma'$ be
the Kneser graphs on the objects of type $J$,~$J'$, respectively,
where $J \subseteq J'$.
Let $(W,R)$ be the Coxeter system of the building, and put
$X = \langle R \setminus J \rangle$,
$X' = \langle R \setminus J' \rangle$, and
$P = BXB$, $P' = BX'B$, so that the vertices of $\Gamma$,~$\Gamma'$
can be viewed as left cosets of $P$,~$P'$, respectively.

\begin{Lemma}
The canonical map $\phi : \Gamma' \to \Gamma$ sending $gP'$ to $gP'P = gP$
is a homomorphism, that is, sends edges to edges.
\end{Lemma}
\Proof
Since $J \subseteq J'$, we have $X' \le X$ and $P' \le P$.
If $gP' \adj hP'$ in $\Gamma'$, that is, if $Bw_0B \subseteq P'g^{-1}hP'$,
then also $Bw_0B \subseteq Pg^{-1}hP$, that is, $gP \adj hP$ in $\Gamma$.
\qed

\begin{Lemma}
Suppose $J^{w_0} = J$. If $a'$ is a vertex of $\Gamma'$ and
$a,b$ are adjacent vertices of $\Gamma$, where $\phi(a') = a$,
then there is a vertex $b'$ of $\Gamma'$, adjacent to $a'$,
with $\phi(b') = b$. That is, edges can be lifted.
\end{Lemma}
\Proof
Let $a' = yP'$, $a = yP$, $b = zP$.
Since $J^{w_0} = J$, we have $w_0X = Xw_0$.
Since $yP,zP$ are adjacent, we have $Bw_0B \subseteq Pz^{-1}yP$,
so that $z^{-1}y \in Pw_0P = BXw_0XB = BXw_0B \subseteq Pw_0B$, and
$w_0B = p^{-1}z^{-1}yB$ for some $p \in P$.
It follows that the vertex $a' = yP'$ is adjacent to $b' = zpP'$.
\qed

\medskip
If the unique coclique extension property fails for objects of some type
$J$ in a building of spherical type with Coxeter system $(W,R)$,
where $J^{w_0} = J$, and $J \subseteq J' \subseteq R$,
then it also fails for objects of type $J'$.

\begin{Proposition}
Let $\Gamma$, $\Gamma'$ be the Kneser graphs on the objects of type
$J$,~$J'$, respectively, in a building of spherical type.
Let $\Sigma$,~$\Sigma'$ be the subgraphs of $\Gamma$,~$\Gamma'$, respectively,
induced on an apartment. If $J^{w_0} = J$, and $J \subseteq J'$,
and the pair $(\Gamma',\Sigma')$ has the unique coclique extension property,
then also $(\Gamma,\Sigma)$ has this property.
\end{Proposition}
\Proof
We can take $\Sigma$,~$\Sigma'$ to be the subgraphs induced on
$WP/P$ and $WP'/P'$.

Let $C$ be a maximal coclique in $\Sigma$. Then $\phi^{-1} C$
is a coclique in $\Sigma'$. If $wP \not\in C$, there is a
$vP \in C$ adjacent to $wP$, and we can lift this edge and find
a neighbour of $wP'$ in $\phi^{-1} C$.
Hence $\phi^{-1} C$ is a maximal coclique in $\Sigma'$.

If $y,z$ are two vertices of $\Gamma$ such that both
$C \cup \{y\}$ and $C \cup \{z\}$ are cocliques, then
both $\phi^{-1} C \cup \phi^{-1}(y)$ and $\phi^{-1} C \cup \phi^{-1}(z)$
are cocliques in $\Gamma'$, and by the unique coclique extension property
of $(\Gamma',\Sigma')$, no vertex of $\phi^{-1}(y)$ is adjacent to
a vertex of $\phi^{-1}(z)$. But then also $y$ and $z$ are nonadjacent.
\qed

\medskip
If the shortest element in $Xw_0X$ equals the shortest element in $X'w_0X'$,
then $\Gamma'$ is a coclique extension of $\Gamma$,
and $(\Gamma',\Sigma')$ has the unique coclique extension property
if and only if $(\Gamma,\Sigma)$ has this property.
For example, if we label the nodes of the $A_n$ diagram by $1, \ldots, n$,
and $J \subseteq \{ 1, \ldots, j\}$ with $\max J = j \le (n+1)/2$, then
$A_{n,J}$ is a coclique extension of $A_{n,j}$ and so has the
unique coclique extension property.

\medskip\noindent
{\tt aeb@cwi.nl}, {\tt jan.draisma@unibe.ch},
{\tt C.Guven@tilburguniversity.edu}


\begin{thebibliography}{99}

\bibitem{bb}
A.~Blokhuis \& A.~E.~Brouwer,
{\it Cocliques in the Kneser graph on line-plane flags in $\PG(4,q)$},
Combinatorica {\bf 37} (2017) 795--804.

\bibitem{BBCFMPS}
A.~Blokhuis, A.~E.~Brouwer, A.~Chowdhury, P.~Frankl, T.~Mussche,
B.~Patk\'os \& T.~Sz\H{o}nyi,
{\it A Hilton-Milner theorem for vector spaces},
Electr. J. Combin. {\bf 17} (2010) R71.

\bibitem{bbg}
A.~Blokhuis, A.~E.~Brouwer \& \c{C}. G\"uven,
{\it Cocliques in the Kneser graph on the point-hyperplane flags
of a projective space},
Combinatorica {\bf 34} (2014) 1--10.

\bibitem{bbs}
A.~Blokhuis, A.~E.~Brouwer \& T.~Sz\H{o}nyi,
{\it On the chromatic number of $q$-Kneser graphs},
Des. Codes Cryptogr. {\bf 65} (2012) 187--197.

\bibitem{bbs2}
A.~Blokhuis, A.~E.~Brouwer \& T.~Sz\H{o}nyi,
{\it Maximal cocliques in the Kneser graph on point-plane flags in $\PG(4,q)$},
Europ. J. Comb. {\bf 35} (2014) 95--104.

\bibitem{bourbaki}
N.~Bourbaki, {\em Lie groups and Lie algebras. Chapters 4--6},
Springer, 2002. 

\bibitem{borel}
A.~Borel, {\em Linear algebraic groups}, Springer, 1991. 

\bibitem{brion-et-al}
M.~Brion \& S.~Kumar, {\it Frobenius Splitting Methods
in Geometry and Representation Theory}, Birkh\"auser, 2005.

\bibitem{qpows}
A.~E.~Brouwer,
{\it The eigenvalues of oppositeness graphs in buildings of spherical type},
pp.~1--10 in: Combinatorics and graphs, Contemp. Math. {\bf 531},
Amer. Math. Soc., Providence, RI, 2010.

\bibitem{chevalley}
C.~Chevalley,
{\it The algebraic theory of spinors},
Columbia University Press, 1954.
Reprinted in: {\it The algebraic theory of spinors and Clifford algebras},
Collected Works Vol. 2, Springer, 1997.

\bibitem{dmm}
J.~De~Beule, S.~Mattheus \& K.~Metsch,
{\it An algebraic approach to Erd\H{o}s-Ko-Rado sets of flags in spherical
buildings},
{\tt arXiv:2007.01104}.


\bibitem{dmw}
J.~D'haeseleer, K.~Metsch \& D.~Werner,
{\it On the chromatic number of two generalized Kneser graphs},
Europ. J. Comb. {\bf 101} (2022) 103474.

\bibitem{EKR}
P.~Erd\H{o}s, C.~Ko \& R.~Rado,
{\it Intersection theorems for systems of finite sets},
Quart. J. Math. Oxford Ser. (2) {\bf 12} (1961) 313--320.

\bibitem{garibaldi-et-al}
S.~Garibaldi, R.~M.~Guralnick \& D.~K.~Nakano,
{\it Globally irreducible Weyl modules},
J. Algebra {\bf 477} (2017) 69--87.

\bibitem{GodsilMeagher}
C.~Godsil \& K.~Meagher,
{\it Erd\H{o}s-Ko-Rado theorems: Algebraic~approaches},
Cambridge Univ. Press, Cambridge, 2015.

\bibitem{GrossNebe2004}
B.~H.~Gross \& G.~Nebe,
{\it Globally maximal arithmetic groups},
J. Algebra {\bf 272} (2004) 625--642.

\bibitem{cicekthesis}
\c{C}.~G\"uven,
{\it Buildings and Kneser graphs},
Ph.~D.~Thesis, Technische Univ. Eindhoven, Eindhoven, 2012.

\bibitem{hiltonmilner}
A.~J.~W.~Hilton \& E.~C.~Milner,
{\it Some intersection theorems for systems of finite sets},
Quart. J. Math. Oxford Ser. (2) {\bf 18} (1967) 369--384.

\bibitem{hsieh}
W.~N.~Hsieh,
{\it Intersection theorems for systems of finite vector spaces},
Discrete Math. {\bf 12} (1975) 1--16.



\bibitem{metsch2019}
K.~Metsch,
{\it An Erd\H{o}s-Ko-Rado result for sets of pairwise non-opposite lines
in finite classical polar spaces},
Forum Math. {\bf 31} (2019) 491--502.

\bibitem{metsch2019a}
K.~Metsch,
{\it An Erd\H{o}s-Ko-Rado result for finite buildings of type $F_4$},
Israel J. Math. {\bf 230} (2019) 813--830.

\bibitem{mw}
K.~Metsch \& D.~Werner,
{\it Maximal cocliques in the Kneser graph on plane-solid flags in $\PG(6,q)$},
Innov. Incidence Geom. {\bf 18} (2020) 39--55.

\bibitem{pepe-et-al}
V.~Pepe, L.~Storme \& F.~Vanhove,
{\it Theorems of Erd\H{o}s-Ko-Rado type in polar spaces},
J. Comb. Th. (A) {\bf 118} (2011) 1291--1312.


\bibitem{Lex}
A.~Schrijver,
{\it Combinatorial Optimization},
Springer, 2003.

\bibitem{Stanton1980}
D.~Stanton,
{\it Some Erd\H{o}s-Ko-Rado theorems for Chevalley groups},
SIAM J. Alg. Disc. Math. {\bf 1} (1980) 160--163.

\bibitem{Tanaka}
H.~Tanaka,
{\it Classification of subsets with minimal width
and dual width in Grassmann, bilinear forms and dual polar graphs},
J. Combin. Th. (A) {\bf 113} (2006) 903--910.

\end{thebibliography}
\end{document}